\documentclass[12pt]{article}

\makeatletter
\newfont{\footsc}{cmcsc10 at 8truept}
\newfont{\footbf}{cmbx10 at 8truept}
\newfont{\footrm}{cmr10 at 10truept}
  \makeatother
\title {\bf 
Bounds on\\
Shannon Capacity
and Ramsey Numbers\\
from Product of Graphs}

\author{
Xiaodong Xu\\
\small Guangxi Academy of Sciences\\[-0.8ex]
\small Nanning, Guangxi 530007, China\\[-0.8ex]
\small \texttt{xxdmaths@sina.com}\\
\\
and\\
\\
Stanis\l aw P. Radziszowski\\
\small Department of Computer Science \\[-0.8ex]
\small Rochester Institute of Technology \\[-0.8ex]
\small Rochester, NY  14623, USA \\[-0.8ex]
\small \texttt{spr@cs.rit.edu}\\\
}

\begin{document}

\def\R{{\cal R}}
\newtheorem{thm}{Theorem}
\newtheorem{col}{Corollary}

\maketitle

\begin{abstract}
In this paper we study Shannon capacity of channels in the context
of classical Ramsey numbers. We overview some of the results
on capacity of noisy channels modelled by graphs, and how some
constructions may contribute to our knowledge of this capacity.

\medskip
We present an improvement to the constructions by Abbott and
Song and thus establish new lower bounds for a special type
of multicolor Ramsey numbers. We prove that our construction
implies that the supremum of the Shannon capacity over all
graphs with independence number 2 cannot be achieved
by any finite graph power. This can be generalized to
graphs with bounded independence number.
\end{abstract}

\bigskip
\bigskip
\bigskip
\noindent
{\bf Keywords:} Shannon channel capacity, Ramsey numbers\\
{\bf AMS classification subjects:} 05C55, 94A24, 05C35


\bigskip
\bigskip
\section{Introduction and Notation}

In this article we study lower bound constructions on some multicolor
Ramsey numbers and their relation to Shannon capacity of noisy channels
modelled by graphs. All
graphs are undirected and loopless, and all colorings are edge-colorings.
The independence number of a graph $G$, i.e. the maximum
number of mutually independent vertices in $G$, will be
denoted by $\alpha(G)$.

\bigskip
For arbitrary graphs $G_1,\ldots,G_n$, where $G_i=(V_i,E_i)$,
we define the {\em graph product}
$G_1 \times \cdots \times G_n$ to be a graph
$G$ on the vertex set $V=V_1 \times \cdots \times V_n$, whose
edges are all pairs of distinct vertices
$\{ (u_1,\ldots,u_n), (v_1,\ldots,v_n)\}$,
such that for each $i$ from 1 to $n$,
$u_i=v_i$ or $\{u_i,v_i\} \in E_i$. This product is
associative, and also commutative up to isomorphisms permuting
the coordinates. $G^n$ will denote the $n$-fold product of the
same graph, namely $G^n=\underbrace{G \times \cdots \times G}_{n}$.
The {\em capacity} $c(G)$ of a graph $G$ was defined by Shannon
\cite{Sha} as the limit

$$c(G) = \lim_{n \rightarrow \infty} \alpha (G^n)^{1/n},\eqno{(1)}$$

\bigskip
\noindent
and is now called the Shannon capacity of a noisy channel modelled by
graph $G$ (see also \cite{Alon1995}, \cite{Alon2006}).
The quantity $c(G)$ is often simply referred to as
the {\em Shannon capacity} of $G$.
The study of $c(G)$ within information theory was initiated by
Shannon \cite{Sha} and has grown to be an extensive
area involving electrical engineering, communication theory,
coding theory, and other fields that
typically use probability theory as a tool.
It may be less known that $c(G)$ attracted
attention of many graph theorists trying to compute it
\cite{Alon2006, Alon2007a, Alon1995, Li, rosta}.
The definitions above, the intuition
below, and our work are representing this graph-theoretic
perspective.

\medskip
Suppose that we have a set $\Sigma$
of $k$ characters which we wish to send over
a noisy channel one at a time. Let $V(G)=\Sigma$, and assume
further that the edges of $G$ indicate a possible confusion
between pairs of characters when transmitted over the channel.
When sending a single character, the maximum number of
characters we can fix, and then choose from for transmission
without danger of confusion, is clearly $\alpha(G)$.
When we use the same channel repeatedly $n$ times, we could
obviously send $\alpha(G)^n$ words of length $n$ by using
an independent set in $G$ at each coordinate. However,
we might be able to do better by sending words from $\Sigma^n$
corresponding to vertices of an independent set of order
$\alpha(G^n)$ in graph $G^n$, in cases when
the general inequality $\alpha(G^n) \ge \alpha(G)^n$ is strict.
The Shannon capacity $c(G)$
measures the efficiency of the best possible strategy when
sending long words over a noisy channel modelled by $G$,
since the limit (1) defining it can be seen as approaching
the effective alphabet size in zero-error transmissions.

\bigskip
A $(k_1,k_2,\ldots,k_n)$-$coloring$, for some $n$ and $k_i \ge 1$,
is an assignment of one of $n$ colors to each edge in a complete graph,
such that the coloring does not contain any monochromatic complete
subgraph $K_{k_i}$
in color $i$, for $1 \le i \le n$. Similarly, a
$(k_1,k_2,\ldots,k_n;s)$-$coloring$
is a $(k_1,\ldots,k_n)$-coloring of the complete graph on $s$ vertices $K_s$.
Let $\R(k_1,\ldots,k_r)$ and $\R(k_1,\ldots,k_n;s)$ denote the set
of all $(k_1,\ldots,k_n)$- and $(k_1,\ldots,k_n;s)$-colorings,
respectively. The Ramsey number $R(k_1,\ldots,k_n)$ is defined to be
the least $s>0$ such that $\R(k_1,\ldots,k_n;s)$ is empty.
In the diagonal case $k_1=\ldots=k_n=k$, we will use simpler
notation $\R_n(k)$ and $\R_n(k;s)$ for sets of colorings and $R_n(k)$
for the Ramsey numbers. The second author maintains a
regularly updated survey \cite{ds1} of the most recent results
on the best known bounds on various types of Ramsey numbers.

\bigskip
In 1971, Erd\H{o}s, McEliece and Taylor \cite{Erdos1971}
were the first to discuss the connections between
$\alpha(G_1 \times \ldots \times G_n)$ and Ramsey numbers. Many papers
followed which studied explicitly Shannon capacity in relation
to independence in product graphs and Ramsey numbers,
like those by Alon et al. \cite{Alon1995,Alon1998,Alon2006,Alon2007b},
Bohman et al. \cite{Boh2003,Boh2009}, and the survey papers
\cite{rosta,Alon2007a,Li}. Here we provide a further link
between lower bounds on some multicolor Ramsey numbers
and Shannon capacity. The result in Theorem 2 of Section 3
enhances our previous constructions from \cite{XXR,XXER}
by establishing new lower bound for a special type of
multicolor Ramsey numbers. This, in turn, implies that
the supremum of the Shannon capacity over all
graphs $G$ with independence number $\alpha(G)=2$
cannot be achieved by using any finite graph power.
The same generalizes to graphs with bounded
independence number.

\bigskip
\section{Some Prior Results}

\bigskip
The main results of a short but interesting paper
by Erd\H{o}s, McEliece and Taylor \cite{Erdos1971}
are summarized in the following theorem.

\bigskip
\noindent
\begin{thm}
{\bf - Erd\H{o}s, McEliece, Taylor - 1971 \cite{Erdos1971}}

\medskip
\noindent
For arbitrary graphs $G_1,\ldots,G_n$,
$$\alpha(G_1 \times \cdots \times G_n)<R(\alpha(G_1)+1,\ldots,\alpha(G_n)+1),\eqno{(2)}$$

\noindent
and for all $k_1,\ldots,k_n>0$ there exist graphs $G_i$ with
$\alpha(G_i)=k_i$, $1\le i \le n$, such that
$$\alpha(G_1 \times \cdots \times G_n)=R(k_1+1,\ldots,k_n+1)-1.\eqno{(3)}$$
\noindent
Furthermore, for the diagonal case $k_i=k$, there exists a single
graph $G$ with $\alpha(G)=k$, such that $\alpha(G^n)=R_n(k+1)-1$.
\end{thm}

\bigskip
This early theorem established
strong links between the Shannon capacity, independence number
of graph products and classical Ramsey numbers. Unfortunately,
all three concepts are notoriously difficult, even for many
very simple graphs. The value of the Shannon capacity of the pentagon,
$c(C_5)=\sqrt{5}$, was computed in a remarkable paper by Lov\'{a}sz
\cite{Lovasz} using tools from linear algebra in a surprising way.
The value of $c(C_7)$ is still unknown, though
significant progress has been obtained by Bohman et al.
for some general cases of odd cycles \cite{Boh2009}
and their complements \cite{Boh2003}.
We only know how to compute $c(G)$ for very special graphs,
like perfect graphs or self-complementary
vertex-transitive graphs, and it seems plausible that even approximating
$c(G)$ may be much harder than {\bf NP}-hard \cite{Alon2006}.

\bigskip
If we use Theorem 1 for non-complete graphs without triangles
in the complement (i.e. $\alpha(G_i)=k_i=2$ for all $i$), then
the Ramsey numbers in question are $R_n(3)$.
It is known that $\lim_{n \rightarrow \infty} R_n(3)^{1/n}$
exists, though it may be infinite. The best established lower
bound for this limit is 3.199... \cite{XXER}. Clearly,
$\lim_{n \rightarrow \infty} (R_n(3)-1)^{1/n} =
\lim_{n \rightarrow \infty} R_n(3)^{1/n}$, and hence by (3) in
Theorem 1, it is equal to the supremum of the Shannon
capacity $c(G)$ over all graphs $G$ with independence
number 2.

Similarly, for any fixed integer $k \ge 3$,
$\lim_{n \rightarrow \infty} R_n(k)^{1/n}$
exists, though again it may be infinite. Furthermore,
this limit is equal to the supremum of the Shannon capacity
$c(G)$ over all graphs $G$ with independence number $k-1$.

\vspace{.2in}
\section{A Ramsey Construction}

\bigskip
This section presents a theorem which gives a new lower bound
construction for some special cases of multicolor Ramsey numbers.
This theorem is improving over an old result by Abbott \cite{Abb}
and Song \cite{Song} that $R_{n+m}(k)>(R_n(k)-1)(R_m(k)-1)$ (see
also \cite{XXER}).
The current approach enhances our previous techniques used in
\cite{XXER,XXR} and summarized in \cite{ds1}.
This result is then linked in Section 4 to the Shannon capacity
of some graphs, in particular graphs with independence 2.

\medskip
We would like to note that a special
product of graphs (and edge-colorings) $G$ and $H$,
denoted $G[H]$, which we used in a few constructions
in \cite{XXER,XXR}, is similar to but distinct from
$G \times H$ usually considered in the context of Shannon capacity.
The vertex set of $G[H]$ is also equal to $V(G) \times V(H)$,
but for graphs $G$ and $H$,
$\{ (u_1,v_1),(u_2,v_2)\}$ is an edge of $G[H]$
if and only if $u_1=u_2$ and $\{ v_1,v_2\} \in E(H)$,
or $\{ u_1,u_2\} \in E(G)$. In the case of colorings,
if $u_1=u_2$ then
$\{ (u_1,v_1),(u_2,v_2)\}$ in $G[H]$
has the same color as $\{ v_1,v_2\}$ in $H$,
else it has the same color as $\{ u_1,u_2\}$ in $G$.
For any edge-coloring $C$, let $C(u,v)$ denote the color
of the edge $\{u,v\}$ in $C$. Thus, equivalently,
for $u_1\not= u_2$ we have
$G[H]((u_1,v_1),(u_2,v_2))=G(u_1,u_2)$,
and $G[H]((u,v_1),(u,v_2))=H(v_1,v_2)$.
Observe that $G[H]$ can be seen as $|V(G)|$ disjoint
copies of $H$ interconnected by many overlapping
copies of $G$. Specifically, there are $|V(H)|^{|V(G)|}$ of them.
Note that, because of this structure, if colors
used in $G$ and $H$ are distinct, then the orders of
the largest monochromatic complete subgraphs in $G[H]$
are the same as in $G$ or $H$, depending on the color.
Finally,
observe that in general the graphs $G[H]$ and
$H[G]$ need not be isomorphic.

\medskip
\bigskip
\begin{thm}
For integers $k, n, m, s \ge 2$, let $G \in \mathcal{R}_n(k;s)$
be a coloring containing an induced subcoloring of $K_m$
using less than $n$ colors. Then
$$
R_{2n}(k) \ge s^2 +m(R_n(k-1,
\underbrace{k, \cdots ,k}_{n-1})-1)+1.
\eqno{(4)}
$$
\end{thm}

\bigskip
\noindent
{\bf Proof. }
Consider coloring $G \in \mathcal{R}_n(k;s)$ with the vertex set
$V(G)=\{ v_1,\ldots ,v_s\}$, and suppose, without loss of generality,
that the set $M=\{ v_1,\dots ,v_m\} , m\le s$, does not induce
any edges of color 1 in $G$. Let $H$ be any critical $n$-coloring
(on the maximum possible number of vertices)
in $\mathcal{R}_n(k-1,k, \cdots ,k)$
with vertices $V(H)=\{ w_1,\dots ,w_t\}$, and hence
$t=\mathcal{R}_n(k-1,k\dots ,k)-1$. In order to prove
the theorem, we will construct a $2n$-coloring
$F \in \mathcal{R}_{2n}(k;s^2+mt)$ with the vertex set
$V(F)=(V(G) \times V(G)) \cup (M \times V(H))$.

\smallskip
We will use colors labeled by integers from 1 to $2n$.
$G$ and $H$ use colors from 1 to $n$, and $F$ from
1 to $2n$.
The structure of coloring $F$ induced on the set
$V(G) \times V(G)$
is similar to that of the special product of $G[G]$,
namely, we set the color of each edge
$e=\{(v_{i_1},v_{i_2}),(v_{j_1},v_{j_2})\}$,
for $1 \le i_1,i_2,j_1,j_2 \le s$, by

\medskip
\bigskip
\[
F(e) = \left\{
\begin{array}{llr}
\mbox{$n+1$} & \mbox{if $i_2=j_2\le m\ $ and $\ G(v_{i_1},v_{j_1})=1$,}\\
\mbox{$G(v_{i_2},v_{j_2})+n$} & \mbox{for $\ i_1 = j_1,$} &
\mbox{\hspace{1.5in} (5)}\\
\mbox{$G(v_{i_1},v_{j_1})$} & \mbox{for other cases with $\ i_1 \not= j_1$.}
\end{array}
\right.
\]

\medskip
\bigskip
\noindent
In addition, the coloring $F$ contains $m$ isomorphic copies of
the coloring $H$ on the vertex sets
$U_i=\{(v_i,w_j)\ |\ 1\le j\le t\}$ for $1\le i \le m$,
each of order $t$. The definition of the coloring of
the edges connecting $U_i$'s follows.

\smallskip
All the edges of the form
$\{(v_{i_1},w_{j_1}),(v_{i_2},w_{j_2})\}$,
for $1\le i_1<i_2\le m$ and $1\le j_1,j_2\le t$,
i.e. the edges between different copies of $H$,
are assigned color $G(v_{i_1},v_{i_2})+n$.
All the edges of the form
$\{(v_{i_1},v_{j_1}),(v_{i_2},w_{j_2})\}$,
for $1\le i_1, j_1 \le s$,
$1\le i_2\le m$, $j_1 \not= i_2$,
and $1 \le j_2 \le t$, are also assigned a high
index color $G(v_{j_1},v_{i_2})+n$.
Finally, the remaining uncolored edges of the form
$\{(v_{i_1},v_q),(v_q,w_{j_2})\}$,
for $1 \le q \le m$, $1 \le i_1 \le s$, and
$1 \le j_2 \le t$,
are assigned color 1.

\bigskip
We will prove that the coloring $F$ constructed
above does not contain any monochromatic $K_k$. 
We already noted that the part of $F$ induced by
the vertices $V(G) \times V(G)$ is similar to $G[G]$.
More precisely, let's denote this part of $F$ by $F'$,
and let $G'$ denote the coloring
obtained from $G$ by renaming all colors from $c$
to $c+n$. Then, if
in $G[G']$ we recolor the edges specified in
the first line of (5) from color 1 to color $n+1$,
then we obtain exactly $F'$.
Next, let the part of $F$ induced by the vertices
$M \times V(H)$ be denoted by $F''$, and the subcoloring
of $G'$ induced by vertices $M$ be denoted by $G''$.
Observe that $F''$ is isomorphic to $G''[H]$.

Since $G[G']$ and $F''$ are both the results of the
special product with different sets of base colors,
they don't contain any monochromatic $K_k$. Furthermore,
since $M$ doesn't induce in $G$ any edges of color 1,
then $F'$ has no monochromatic $K_k$ either.
Thus, if there is
a monochromatic $K_k$ in $F$ it must intersect both
$V(G) \times V(G)$ and $M \times V(H)$.
Next, it is not hard
to see that the structure of $F'$ and $F''$ and how
they swap the roles of colors with labels
$\le n$ and $>n$ prevent any monochromatic
$K_k$ in all colors different from 1 and $(n+1)$.
Now, note that the lastly added edges in color 1
between $F'$ and $F''$ join the blocks of vertices
with the same index $v_q$ in one position, and
one can conclude that no monochromatic $K_k$
in color 1 can arise because
$H \in \mathcal{R}_n(k-1,k, \cdots ,k)$.
Finally, no $K_k$ is
formed in color $(n+1)$ since $G''$
has no edges in color $(n+1)$.
This completes the proof.
$\diamondsuit$

\vspace{.3in}
We wish to comment that Theorem 2 with any lower bound better
than $s^2$ would be sufficient for the results in Section 4.
Observe that the required subcoloring with $m>0$ exists
in all nontrivial cases.

\bigskip
\section{Shannon Capacity}

\bigskip
It can be easily shown that $R_{2n}(3) > (R_n(3)-1)^2$,
for example by using inequalities (7) or (12) in \cite{XXER}.
Now, this can be improved by the construction of Theorem 2,
as in the corollary below. This corollary is interesting in itself
since it improves over the previous lower bound recurrence
on $R_n(3)$, but first of all it will let us prove Theorem 3
on Shannon capacity of graphs with independence number 2.

\bigskip
\noindent
{\bf Corollary}
{\em For all integers $n \ge 2$,}

\medskip
$R_{2n}(3) \ge (R_n(3)-1)^2 + m(R_{n-1}(3)-1)+1$, for
$m = \lceil(R_n(3)-2)/n\rceil$.

\bigskip
\noindent
{\bf Proof. }
Each vertex of any coloring in $\R_n(3;s)$ has at least
$m = \lceil (s-1)/n \rceil$ neighbors in the same color, which must
induce a coloring in $\R_{n-1}(3;m)$.
Theorem 2 with $k=3$ implies the claim.
$\diamondsuit$

\bigskip
\begin{thm}
If the supremum of the Shannon capacity $c(G)$
over all graphs with independence number $2$ is finite
and equal to $C$, then $C>\alpha (G^n)^{1/n}$
for any graph $G$ with
independence number $2$ and for any positive integer $n$.
\end{thm}

\noindent
{\bf Proof. }
Suppose that
$C$ is achieved by some graph $G$ with $\alpha(G)=2$,
and hence by (2) we have $C^n=\alpha(G^n)<R_n(3)$.
By the second part of Theorem 1, we know that there exists
a graph $H$ with independence number $2$ such that
$\alpha(H^{2n})=R_{2n}(3)-1$, and by Theorem
2 we see that $\alpha(H^{2n})>(R_n(3)-1)^2$.
This contradicts the fundamental inequality
$\alpha(G_1 \times G_2)\ge\alpha(G_1)\alpha(G_2)$,
Theorem 1 and the assumption that $C$ is realized by $G$.
$\diamondsuit$

\bigskip
Observe that in the case of infinite supremum,
$\lim_{n \rightarrow \infty} R_n(3)^{1/n}$ must also be infinite.
In other words, together with Theorem 3, this means that
the supremum of the Shannon capacity over all
graphs $G$ with independence number $\alpha(G)=2$
cannot be achieved by any finite graph power.

\bigskip
It is not difficult to generalize Theorem 3 for $k\ge 3$
to $\alpha(G)=k-1$, $R_n(k)$ and
the supremum of the Shannon capacity over all graphs with
independence number $k-1$, as stated in the following
Theorem 4. We omit the details which are
analogous to those in Corollary and Theorem 3.

\bigskip
\begin{thm}
$(a)$
For every positive integer $n_0$,
$R_{n_0}(k)^{1/{n_0}} < lim _ {n\rightarrow \infty} R_n(k)^{1/n}$, and\break
$(b)$
the supremum of the Shannon capacity over all graphs with
bounded independence number cannot be achieved by any
finite graph power.\\
\end{thm}

\bigskip
\section*{Acknowledgments}

The work of the first author was partially supported by
the Guangxi Natural Science Foundation (2011GXNSFA018142).
We would like to thank Alexander Lange for his careful
reading and suggestions of improvements,
and finally we are very grateful to the reviewer whose
insightful comments led to a much better presentation
of this paper.


\end{document}